\documentclass[a4paper,reqno]{amsart}
\usepackage{amssymb}
\usepackage{latexsym}
\usepackage{amsmath}
\usepackage{euscript}
\usepackage{graphics,color}
\usepackage[all]{xy}
\usepackage[margin=3cm]{geometry}
\usepackage{MnSymbol} % for \udots

\newcommand{\be}{\begin{equation}}
\newcommand{\ee}{\end{equation}}
\newcommand{\ba}{\begin{eqnarray}}
\newcommand{\ea}{\end{eqnarray}}
\newcommand{\baa}{\begin{eqnarray*}}
\newcommand{\eaa}{\end{eqnarray*}}
\newcommand{\bb}{\begin{thebibliography}}
\newcommand{\eb}{\end{thebibliography}}

\newcommand{\bi}[1]{\bibitem{#1}}
\newcommand{\lab}[1]{\label{#1}}
\newcommand{\re}[1]{(\ref{#1})}

\newcounter{my}
\newcommand{\he}%
   {\stepcounter{equation}\setcounter{my}%
   {\value{equation}}\setcounter{equation}0%
   }%
\newcommand{\she}%
   {\setcounter{equation}{\value{my}}%
    }%

\renewcommand\t{\tilde}

\newcommand{\olsi}[1]{\,\overline{\!{#1}}} % overline short italic

\newtheorem{pr}{Proposition}

\theoremstyle{definition}

\numberwithin{equation}{section}

\begin{document}

\title[circle Jacobi algebra]{The CMV bispectrality of the Jacobi polynomials on the unit circle}

\author{Luc Vinet}

\address{IVADO and Centre de Recherches Math\'ematiques, Universit\'e de Montr\'eal, P.O. Box 6128, Centre-ville Station,
Montr\'eal (Qu\'ebec), H3C 3J7, Canada}

\author{Alexei Zhedanov}

\address{School of Mathematics, Renmin University of China, Beijing 100872, China}

%\vspace*{5mm}

\begin{abstract}
We show that the Jacobi polynomials that are orthogonal on the unit circle (the Jacobi OPUC) are CMV bispectral. This means that the corresponding Laurent polynomials in the CMV basis satisfy two dual ordinary eigenvalue problems: a recurrence relation and a differential equation of Dunkl type. This is presumably the first nontrivial explicit example of CMV bispectral OPUC. We introduce the circle Jacobi algebra which plays the role of hidden symmetry algebra for the Jacobi OPUC. All fundamental properties of the Jacobi OPUC can be derived from representations of this algebra.
\end{abstract}

\maketitle
To the memory of André Ronveaux who had a great influence on the field of OPSFA and its community and who maintained friendly ties with Montreal.

\section{Introduction}
\setcounter{equation}{0}

In this paper we show that the polynomials of Jacobi type that are orthogonal on the unit circle (the Jacobi OPUC) are ``classical" in the sense that this adjective has for the orthogonal polynomials of the Askey scheme. Introduced by Szeg\H{o} \cite{Szego} roughly 100 years ago as ``circle" analogs of the ordinary Jacobi polynomials, the Jacobi OPUC $\Phi_n(z)$ have several properties that make them comparable to their real line analogs \cite{Badkov}; their orthogonality relation on the unit circle is fully spelled out, their recurrence (Verblunsky) coefficients $a_n$ are completely known, and they have an explicit expression in terms of the Gauss hypergeometric function \cite{DVZ}. However, one point in this parallel has been missing, and boils down to the question of whether or not the Jacobi OPUC are bispectral or more specifically ``classical".

We shall here finally settle this issue in the affirmative and identify the bispectrality property of the Jacobi OPUC. 
We shall show in this respect that these polynomials satisfy two dual eigenvalue problems: one where the operator acts on the degree variable $n$ (the recurrence CMV relation) and the other involving an operator affecting the argument $z$ (the differential eigenvalue equation). Moreover, we shall introduce the Jacobi OPUC algebra which will be seen to play the role of ``hidden symmetry" algebra for the Jacobi OPUC. Finally, we shall indicate how the circle Jacobi algebra is related to the ordinary Jacobi algebra \cite{GIVZ}.

The paper is organized as follows. In Section 2 we recall basic facts from the theory of orthogonal polynomials on the unit circle (OPUC). In Section 3 we review the CMV approach to OPUC and introduce the reflection operators ${\mathcal M}_1, {\mathcal M}_2$ that will play an  important role in our study. Section 4 offers a brief description of the Szeg\H{o} mapping from OPUC to orthogonal polynomials on the real line (PRL). The bispectral problem for OPUC in the CMV basis is formulated in Section 5 and compared with the well known bispectrality property of the ``classical" PRL. In section 6, as a lead towards our central result, we consider the simple example of the single-moment OPUC which have an elementary explicit description. The main observation of this section is that these polynomials are CMV-bispectral. We also note in this context the appearance of (a special case of) an algebra to be called  the ``circle Jacobi algebra". We study this algebra in its general form in Section 7 and, assuming that the associated OPUC are CMV bispectral, we derive the expression of their Verblunsky parameters $a_n$ from the sole use of the defining relations of this algebra. We find that the parameters $a_n$ thus obtained correspond to those of the Jacobi OPUC. In Section 8, we introduce a first order differential operator $\mathcal K$ of Dunkl type to see that it has the Jacobi OPUC as eigenfunctions. We further show that this operator together with ${\mathcal M}_1, {\mathcal M}_2$ provides a realization of  the circle Jacobi algebra.  The relation with the non-symmetric Jacobi polynomials studied in \cite{KB} will also be pointed out. In Section 9 we demonstrate that a central extension of the (quadratic) Jacobi algebra is embedded in the circle Jacobi algebra and explain how the operator $\mathcal K$ is related to the Gauss hypergeometric operator (which acts diagonally on the Jacobi PRL). The last section is a brief conclusion.

\section{Orthogonal polynomials on the unit circle - OPUC}
\setcounter{equation}{0}

With a positive measure $d \mu(\theta)$ on the unit circle, one can define the trigonometric moments 
\be
\sigma_n = \frac{1}{2 \pi} \, \int_{0}^{2 \pi} e^{i \theta n} d \mu(\theta) , \quad n=0, \pm 1, \pm 2, \dots  \lab{g_n_def} \ee
As an example, for the simple Lebesgue measure $d \mu(\theta) = d \theta$, all moments vanish apart from $\sigma_0$
\be
\sigma_0= 1, \quad \sigma_n =0, \; n=\pm 1, \pm 2, \dots \lab{Leb_g} \ee
Using the trigonometric moments $\sigma_n$, one can construct the Toeplitz determinants \be \Delta_n= \left |
\begin{array}{cccc} \sigma_0 & \sigma_{1} & \dots &
\sigma_{n-1}\\ \sigma_{-1}& \sigma_0 & \dots & \sigma_{n-2}\\ \dots & \dots & \dots & \dots\\
\sigma_{1-n} & \sigma_{2-n} & \dots & \sigma_0 \end{array} \right
| , \quad n=1,2,\dots \lab{Delta} \ee 
The positivity of the measure $d \mu(\theta)$ implies the obvious property
\be
\sigma_{-n} = \bar{\sigma}_n,
\lab{sym_sigma} \ee
where $\bar{\sigma}_n$ stands for the complex conjugate of ${\sigma}_n$.
It is well known that the positivity of the measure $d \mu(\theta)$ is equivalent to the positivity of  Toeplitz determinants \cite{Simon}
\be
\Delta_n >0, \quad n=1,2,3, \dots \lab{pos_Delta} \ee

Given a measure $\mu(\theta)$, one can introduce the polynomials \cite{Simon}
\ba \nonumber
&&\Phi_n(z)=(\Delta_n)^{-1} \left |
\begin{array}{cccc} \sigma_0 & \sigma_1 & \dots & \sigma_n \\ \sigma_{-1} & \sigma_0 &
\dots & \sigma_{n-1} \\ \dots & \dots & \dots & \dots\\
\sigma_{1-n}& \sigma_{2-n}& \dots & \sigma_1\\ 1& z & \dots & z^n
\end{array} \right |. \lab{deterPhi} \ea
It is seen that the polynomial $\Phi_n(z)$ is a monic polynomial of degree $n$:
\be
\Phi_n(z) = z^n + O(z^{n-1}). \lab{Phi_monic} \ee
By construction, the polynomials $\Phi_n(z)$ obey the orthogonality property
\be
 \int_{0}^{2 \pi} \Phi_n(e^{i \theta}) e^{-i j \theta } d \mu(\theta) = 0, \; j=0,1,\dots, n-1,
 \lab{ort_1} \ee
which is equivalent to \cite{Simon}
\be
\int_{0}^{2 \pi} \Phi_n(e^{i \theta}) \olsi{\Phi_m(e^{i \theta})} d \mu(\theta) = \int_{0}^{2\pi} \Phi_n(e^{i \theta}) \olsi{\Phi}_m(e^{-i \theta}) d \mu(\theta)= h_n \, \delta_{nm}, \lab{ort_2} \ee
where
\be
h_n = \frac{\Delta_{n+1}}{\Delta_n}>0. \lab{h_Delta} \ee
The notation $\olsi{\Phi}_m(z)$ represents the polynomial that is obtained by taking the complex conjugate of the coefficients of  $\Phi_n(z)$ and leaving the argument $z$ unchanged. Polynomials $\Phi_n(z)$ defined as above and verifying \eqref{ort_1} are called orthogonal polynomials on the unit circle, or OPUC for brevity.
They satisfy the fundamental Szeg\H{o} recurrence relation \cite{Simon}
\be
\Phi_{n+1}(z) = z \Phi_n(z) - \bar a_n \Phi_n^*(z), \lab{Sz_rec} \ee 
where
\be
\Phi_n^*(z) = z^n \olsi{\Phi}_n(z^{-1}).
\lab{Phi*} \ee
The Verblunsky parameters $a_n$ satisfy the condition
\be
|a_n| <1 . \lab{a<1} \ee
One can show \cite{Simon} that the recurrence relation \re{Sz_rec} together with the condition \re{a<1} is necessary and sufficient for the polynomials $\Phi_n(z)$ to be OPUC with respect to a positive measure $d \mu(\theta)$. The normalization constants $h_n$ have the following expression in terms of the parameters $a_n$ \cite{Simon}:
\be
h_0=1, \; h_n = \left(1-|a_0|^2 \right) \left(1-|a_1|^2 \right) \dots \left(1-|a_{n-1}|^2 \right), \quad n=1,2,\dots
\lab{h_a} \ee

%\vspace{5mm}
An important special case is when the parameters $a_n$ are {\it real} and satisfy $-1<a_n<1$. In this case the orthogonality relation reads
\be
\int_{-\pi}^{\pi} \Phi_n(\exp(i \theta)  \Phi_m(\exp(-i \theta) \rho(\theta) d \theta = h_n \delta_{nm}
\lab{ort_real} \ee
with $\rho$ a positive and symmetric weight function:
\be
\rho(-\theta) = \rho(\theta).
\lab{rho_sym} \ee
The {\it orthonormal} OPUC $\varphi_n(z)$ are defined as
\be
 \Phi_n(z) = \sqrt{h_n}  \varphi_n(z).
 \lab{def_vphi} \ee
They satisfy the orthogonality relation
\be
\int_{-\pi}^{\pi} \varphi_n(\exp(i \theta))  \varphi_m(\exp(-i \theta)) \rho(\theta) d \theta = \delta_{nm}.
\lab{ort_varphi} \ee

\section{OPUC and the CMV basis}
\setcounter{equation}{0}
Introduce the following CMV basis of the space of Laurent polynomials
\be 
\chi_0 =1, \: \chi_1(z) =z, \: \chi_2(z) = z^{-1}, \dots, \chi_{2n-1}(z) = z^{n}, \: \chi_{2n}=z^{-n}, \dots
\lab{CMV} \ee 
The dual CMV basis is defined by the involution
\be
\chi^{\star}_n(z) = \chi_n(1/z).
\lab{chi_star} \ee
In what follows we will restrict to the case of real $a_n$. 
The Laurent polynomials $\psi_n(z)$ are defined as 
\be
\psi_{2n}(z) = z^n \Phi_{2n}(1/z), \; \psi_{2n+1}(z) = z^{-n} \Phi_{2n+1}(z).
\lab{psi_def} \ee
Clearly, for any $n=0,1,2,\dots$ the Laurent polynomial  $\psi_n(z)$ is a linear combination of $\chi_0, \chi_1(z), \\
\dotsm ,\chi_n(z)$. These $\psi_n(z)$ satisfy specific recurrence relations of CMV type \cite{Simon}. In order to derive them let us introduce the pair of reflection operators $R$ and $ZR$ defined as
\be
Rf(z) = f(1/z), \quad ZR(f)(z)= zf(1/z), \quad Z f(z) = zf(z).
\lab{RZR} \ee
Obviously
\be
R^2=(ZR)^2=\mathcal{I}.
\lab{M_12_I} \ee

It is straightforwardly verified that the action of the operator $R$ on the basis $\psi_n$ can be presented as
\be
R {\vec \psi(z)} = {\vec \psi(1/z)} = {\mathcal M}_1 {\vec \psi(z)}
\lab{RM1} \ee
where $\vec \psi(z)$ is the infinite-dimensional vector ${\vec \psi(z)}=\{\psi_0(z), \psi_1(z), \dots, \}$ and where $\mathcal{M}_1$ is the semi-infinite block-diagonal matrix
\be
\mathcal{M}_1 =
 \begin{pmatrix}
1 \\
& a_{1} & 1 &  &    \\
  &1-a_1^2 & -a_{1} &  &   \\
   & &  &               a_{3} & 1  \\
   & & &              1-a_3^2 & -a_{3}  \\
  &  & &   & &           a_{5} & 1  \\
   &  & &   & &         1-a_5^2 & - a_{5}  \\
&   &  & & & & & \ddots  \\
 \end{pmatrix}.
\lab{M1_def} \ee 
Similarly
\be
ZR {\vec \psi(z)} = {z \vec\psi(1/z)} = {\mathcal M}_2 {\vec \psi(z)}
\lab{RM2} \ee
with ${\mathcal M}_2$ the block-diagonal matrix
\be
 \mathcal{M}_2 =
 \begin{pmatrix}
  a_{0} & 1 &  &    \\
  1-a_0^2 & -a_{0} &  &   \\
   &  &              a_{2} & 1  \\
   &  &              1-a_2^2 & -a_{2}  \\
  &  &    & &          a_{4} & 1  \\
   &  &    & &         1-a_4^2 & - a_{4}  \\
&   &  & & & & \ddots  \\
 \end{pmatrix}.
\lab{M2_def} \ee
It follows that the vector ${\vec \psi}$ satisfies the generalized eigenvalue problem 
\be
\mathcal{M}_2 { \vec  \psi(z)} = z \mathcal{M}_1 { \vec  \psi(z)}.
\lab{GEVP_psi} \ee
Multiplying \re{GEVP_psi} by $\mathcal{M}_1$, we find that $\vec \psi(z)$ satisfies also the {\it ordinary} eigenvalue problem 
\be 
{\mathcal C} \vec \psi = z \vec \psi
\lab{CMV_rec} \ee
where $\mathcal C = \mathcal{M}_1 \mathcal{M}_2$ is a pentadiagonal (CMV) matrix \cite{Simon}.

\section{Szeg\H{o} pair of orthogonal polynomials on the real line}
\setcounter{equation}{0}
Szeg\H{o} introduced \cite{Szego} a map from the polynomials $\Phi_n(z)$ on the unit circle with real coefficients $a_n$ to polynomials $P_n(x(z))$ orthogonal on the interval $[-2,2]$ of the real line. Explicitly, the Szeg\H{o} map is given by the formulas
\be
P_0(x(z)) =1, \quad P_n(x(z)) = z^{1-n} \Phi_{2n-1}(z) + z^{n-1} \Phi_{2n-1}(1/z), \; n=1,2,\dots ,
\lab{P_Phi} \ee   
where 
\be
x(z) = z+1/z.
\lab{x_z} \ee
It is seen that $P_n(x(z))$ are monic polynomials of the argument $x(z)$ which satisfy the recurrence relation
\be
P_{n+1}(x) + b_n P_{n}(x) + u_n P_{n-1}(x) = xP_n(x) 
\lab{rec_P_Sz} \ee
where the recurrence coefficients are
\ba
u_n = (1+a_{2n-1}) (1-a_{2n-3}) (1-a_{2n-2}^2), \; b_n = a_{2n}(1-a_{2n-1})-a_{2n-2}(1+a_{2n-1}).
\lab{xi_Sz} \ea
The condition $a_{-1}=-1$ is assumed in these formulas. 

Moreover it is convenient to introduce the ``companion" polynomials $Q_n(x(z))$ \cite{Szego} that are orthogonal on the same interval $[-2,2]$: 
\be
Q_n(z) = \frac{z^{-n} \Phi_{2n+1}(z) - z^{n} \Phi_{2n+1}(1/z)}{z-z^{-1}},
\lab{Q_Sz} \ee
and which satisfy the recurrence relation 
\be
Q_{n+1}(x) + \t b_n Q_{n}(x) + \t u_n Q_{n-1}(x) = xQ_n(x) 
\lab{rec_Q_Sz} \ee
with
\ba
\t u_n = (1+a_{2n-1})(1-a_{2n+1})(1-a_{2n}^2), \; \t b_n = a_{2n}(1-a_{2n+1}) -a_{2n+2}(1+a_{2n+1}).
\lab{xi_Sz_Q} \ea

The polynomials $Q_n(x(z))$ can be obtained from $P_n(x(z))$ by the double Christoffel transform at the two spectral points $x= \pm 2$. This leads to the following formula expressing $Q_n(x(z))$ in terms of $P_n(x(z))$:
\ba
&&\left(z-z^{-1}\right)^2 Q_{n-1}(x(z)) = P_{n+1}(x(z)) + (a_{2n}+a_{2n-2})(1-a_{2n-1})P_n(x(z)) \nonumber \\ &&-(1-a_{2n-1})(1-a_{2n-3})\left(1-a_{2n-2}^2 \right)P_{n-1}(x(z)). \lab{Q->P} \ea
Using the recurrence relation \re{rec_P_Sz}, one can present \re{Q->P} in the equivalent form
\ba
\left(z-z^{-1}\right)^2 Q_{n-1}(x(z)) = \left( z+z^{-1} + 2a_{2n-2} \right) P_n(x(z)) - 2\left(1-a_{2n-3} \right) \left( 1-a^2_{2n-2}\right)P_{n-1}(x(z)). \lab{Q-->P} \ea
The mapping from $Q_n(x)$ to $P_n(x)$ is given by the Geronimus transformation
\be
P_n(x) = Q_n(x) -(1+a_{2n-1})(a_{2n}+a_{2n-2})Q_{n-1}(x) -(1+a_{2n-1})(1+a_{2n-3})\left( 1-a_{2n-2}^2 \right)Q_{n-2}(x).
\lab{PQG} \ee
The inverse map from the polynomials $P_n(x), Q_n(x)$ to $\Phi_n(z)$ reads
\ba
&&\psi_{2n-1}(z) = \frac{1}{2} \left[ P_n(x(z)) +\left( z-z^{-1} \right) Q_{n-1}(x(z)) \right] , n=1,2,\dots \nonumber \\
&&\psi_{2n}(z) = \frac{1}{2} \left[ (1-a_{2n-1}) P_n(x(z)) -(1+a_{2n-1})\left( z-z^{-1} \right) Q_{n-1}(x(z))\right], \; n=0,1,2,\dots 
\lab{psi_PQ}
\ea
or equivalently, 
\ba
&&\psi_{2n-1}(z) = \frac{(z+a_{2n-2})P_n(x(z)) - (1-a_{2n-3})(1-a^2_{2n-2})P_{n-1}(x(z))}{z-z^{-1}}   , \: n=1,2,\dots \nonumber \\
&&\psi_{2n}(z) = \frac{\left(a_{2n-1}z +z^{-1} + a_{2n-2}(1+a_{2n-1}) \right)P_n(x(z))}{z-z^{-1}} - \nonumber \\ 
&&-\frac{(1+a_{2n-1})(1-a_{2n-3})(1-a^2_{2n-2})P_{n-1}(x(z))}{z-z^{-1}}, \; n=0,1,2,\dots 
\lab{psi_PP}
\ea
From these formulas we have another useful expression of the polynomials $P_n(x)$ and $Q_n(x)$ in terms of $\psi_n(z)$:
\ba
&&P_n(x(z)) = \psi_{2n}(z) +\left( 1+a_{2n-1}\right) \psi_{2n-1}(z), \nonumber \\
&&\left(z-z^{-1}\right)Q_{n-1}(x(z)) = -\psi_{2n}(z) +\left( 1- a_{2n-1}\right) \psi_{2n-1}(z).
\lab{PQ_psi} \ea

\section{Bispectrality}
\setcounter{equation}{0}
One is presented with a bispectral problem when a function $f(x,y)$ in two variables $x,y$ satisfies simultaneously two different eigenvalue equations:
\be
V_x f(x,y) = A(y) f(x,y), \quad U_y f(x,y) = B(x) f(x,y),
\lab{VU_bispec} \ee
where $V_x$ is a differential or difference operator acting on the variable $x$ while $U_y$ is an operator acting on the variable $y$. $A(y)$ and $B(x)$ are the corresponding eigenvalues. A classic paper on the topic is \cite{DG}. Such problems arise in various contexts and one of special interest here is provided by orthogonal polynomials $P_n(x)$ on the real line. Let us elaborate a little. 

The orthogonal PRL are completely characterized by their three-term recurrence relation
\be
P_{n+1}(x) + b_n P_n(x) + u_n P_{n-1}(x) = xP_n(x)
\lab{rec_P} \ee
with $b_n,u_n$ real coefficients and $u_n>0, \: n=1,2,\dots$.
Taking for initial conditions 
\be
P_{0}(x)=1, \quad P_{-1}=0,
\lab{ini_P} \ee
the recurrence relation \re{rec_P}  determines a unique system of monic polynomials $P_n(x) = x^n + O(x^{n-1})$ that are orthogonal on the real line 
\be
\int_{\alpha}^{\beta} P_n(x) P_m(x) d \mu(x) = h_n \delta_{nm}, \quad h_n = u_1 u_2 \dots u_n >0.
\lab{ort_P} \ee
The recurrence relation \re{rec_P} can be presented in algebraic form as the eigenvalue problem
\be
J {\vec P} = x \vec P,
\lab{JPP} \ee
where $\vec P = \{P_0, P_1(x), \dots, P_n(x), \dots, \}$ is generically a vector in an infinite-dimensional Hilbert space, and where $J$ is a Jacobi (i.e. tri-diagonal) matrix with entries $b_0,b_1, \dots$ on the main diagonal, $u_1,u_2, \dots$ on the bottom sub-diagonal and 1s everywhere on the upper sub-diagonal.

A family of orthogonal polynomials will be {\it bispectral} if apart from the obeying the eigenvalue equation \re{JPP} that is imposed by their orthogonality, the polynomials $P_n(x)$ satisfy in addition a ``dual" eigenvalue problem
\be
W P_n(x) = \lambda_n P_n(x),
\lab{WPP} \ee
where $W$ is either a differential or difference linear operator acting on the argument $x$.

If the operator $W$ in \eqref{WPP} is of second order, the bispectral polynomials $P_n(x)$ are furthermore said to be {\it classical}. The exhaustive list of classical polynomials that are eigenfunctions of a second order differential operators consists of the Jacobi, Laguerre, Hermite and Bessel polynomials.\footnote{Outside the realm of classical polynomials, the (bispectral) Krall polynomials obey equations of the form \eqref{WPP} with $W$ a differential operator of degree $N>2$. The classification of these OPs has been undertaken: it has been completed for $N=4$ while for $N>4$ apart from conjectures the problem remains open.} The complete set of the classical polynomials that verify \eqref{WPP} with $W$ is a difference second order operator defined on some grid is also known and made out of the Askey-Wilson polynomials and their special and degenerate cases (e.g. the Racah, Hahn, Meixner,... polynomials), see e.g. \cite{VZ_Bochner}. In summary, all the hypergeometric orthogonal polynomials forming the Askey tableau are bispectral and in fact classical.

It has moreover been shown that the fundamental properties (recurrence relation, difference equation, duality etc) of all these families of classical polynomials \cite{KLS} can be derived from the representations of the Askey-Wilson algebra \cite{Zhe_AW} and its degenerate forms and special cases (Racah, Bannai-Ito, Hahn, Jacobi). See for illustration the reviews \cite{GVZ}, \cite{DGTVZ}.

In the case of OPUC $\Phi_n(z)$, the notion of bispectrality seems a priori to require an extension to generalized eigenvalue problems. Indeed, the recurrence relation 
\be
\Phi_{n+1}(z) + \xi_n \Phi_n(z) = z\left( \Phi_n(z) + \eta_n \Phi_{n-1}(z)\right) 
\lab{rec_Phi} \ee
that the polynomials $\Phi_n(z)$ satisfy can be presented as the generalized eigenvalue equation
\be
J_1 {\vec \Phi} = z J_2 {\vec \Phi}
\lab{GEVP_Phi} \ee
for some bi-diagonal matrices $J_1$ and $J_2$. One could then say that the OPUC $\Phi_n(z)$ are bispectral if they obey a ``dual" eigenvalue equation of the same generalized type as \eqref{GEVP_Phi}, i.e. 
\be
W_1 \Phi_n(z) = \lambda_n W_2 \Phi_n(z),
\lab{WW_Phi} \ee
where $W_1$ and $W_2$ are some differential or difference operators. Some preliminary results in this direction were recently obtained in \cite{VZ_Askey} but the general classification of OPUC that are bispectral in that sense is for now proving more complicated than in the case of OPRL.

We can however take a different viewpoint by not formulating the bispectrality of OPUC through equations satisfied by the polynomials $\Phi_n(z)$ themselves but by focusing instead on the Laurent polynomials $\psi_n(z)$ defined by \re{psi_def}.
The advantage of such an approach is that the Laurent polynomials $\psi_n(z)$ satisfy the {\it ordinary} eigenvalue problem \re{CMV_rec} with some pentadiagnal matrix $\mathcal C$. Therefore, in situations where the Laurent polynomials $\psi_n(z)$ obey in addition a dual standard spectral problem
\be
W \psi_n(z) = \lambda_n \psi_n(z), \; n=0,1,2,\dots
\lab{W_psi} \ee
with $W$ a differential or difference operator acting on the variable $z$, we will say that the OPUC $\Phi_n(z)$ exhibit bispectrality of CMV type.
Gr\"unbaum and Vel\'asquez considered in \cite{GV} the case when $W$ is a linear differential operator of arbitrary order. Their main result is that the only OPUC with such CMV bispectrality are the trivial ones: $\Phi_n(z) = z^n$ (that Simon calls the ``free case" OPUC). In this instance, the operator arising in \re{W_psi} is simply $W=z\partial_z$.

We shall show in the following that nontrival examples of CMV bispectral OPUC can be found if one allows $W$ to be Dunkl type operators, that is differential-difference operators containing {\it reflections}. We shall focus on the Jacobi OPUC. The corresponding Dunkl operator will be identified by positing through an educated guess a three-generated algebra whose representations will be shown to provide a full characterization of the Jacobi OPUC including the explicit expression of the Verblunsky parameters $a_n$, the operator $W$ and its eigenvalues $\lambda_n$

%\newpage

\section{An informative special case: the single moment OPUC}
\setcounter{equation}{0}
Consider the OPUC with  the weight function
\be
\rho(\theta) = \frac{1}{2 \pi} \left(1-\xi \cos \theta \right)
\lab{rho_SM} \ee
where $\xi$ is a real parameter.
Only one moment (apart from $\sigma_0=1$) is nonzero:  $\sigma_1= \sigma_{-1}=-\xi/2$ and hence the name single moment OPUC \cite{Simon}.
In what follows we shall consider the simplest situation where $\xi=1$. In this case it is easily verified that
\be
a_n =-\frac{1}{n+2} 
\lab{aSM} \ee
and that the OPUC have the explicit expression \cite{Simon}:
\be
\Phi_n(z) = \frac{1}{n+1} \: \sum_{k=0}^n (k+1)z^k.
\lab{Phi_SM} \ee
Introduce the Dunkl type operator
\be
\mathcal{K} = z \partial_z + \frac{z}{1-z} \left(R-\mathcal{I} \right).
\lab{L_SM} \ee
It is immediately checked that 
\be
\mathcal{K} \psi_n(z) = \lambda_n \psi_n(z),
\lab{EIG_SM} \ee
with
\be
\lambda_n = \left\{ {-n/2, \quad \mbox{if} \; n \; \mbox{even}  \atop (n+3)/2, \quad \mbox{if} \; n \; \mbox{odd}}.  \right.
\lab{lambda_SM} \ee
Moreover, it can be easily verified that the operator $L$ is self-adjoint
\be
\rho^{-1}(z) \mathcal{K}^{\dagger} \rho(z) ={\mathcal K}.
\lab{SA_SM} \ee

Now consider the three operators: $\mathcal{K}, {\mathcal M}_1 = R, {\mathcal M}_2 = ZR$. It is easy to see that these operators verify the algebraic relations: 
\be
\{{\mathcal K}, {\mathcal M}_1 \} = {\mathcal M}_1 - \mathcal{I},
\lab{SM_1A} \ee
and 
\be
\{{\mathcal K}, {\mathcal M}_2 \} = 2{\mathcal M}_2 + \mathcal{I},
\lab{SM_2A} \ee
where $\{X,Y\}=XY+YX$ denotes the anticommutator.
Let us mention in passing, that the algebraic relations \re{SM_1A}-\re{SM_2A} can be checked in the two ``dual" pictures by either using in the argument representation, ${\mathcal K}$ given by the Dunkl type differential operator \eqref{L_SM} and ${\mathcal M}_1$ and ${\mathcal M}_2$ by the reflection idempotents acting on the space of functions $f(z)$ in $z$ or, in the degree representation, taking for ${\mathcal K}$ the diagonal matrix $diag(\lambda_0, \lambda_1, \dots, \lambda_n, \dots)$ made out of the eigenvalues and for  ${\mathcal M}_1$ and ${\mathcal M}_2$ the block-diagonal matrices \re{M1_def}-\re{M2_def}.

\section{The circle Jacobi algebra and representations}
\setcounter{equation}{0}
In the previous section we have demonstrated that the single moment OPUC are bispectral and that the associated operators $\mathcal{K}, {\mathcal M}_1, {\mathcal M}_2$ (recall that $\mathcal{C} = {\mathcal M}_1 {\mathcal M}_2$) satisfy a pair of simple algebraic relations \re{SM_1A}-\re{SM_2A}. This brings us to consider a natural generalization of the algebra thus defined. 

Keeping the same notation $\mathcal{K}, {\mathcal M}_1, {\mathcal M}_2$ for its three generators, this new ``algebra" is defined by the relations
\be
\{ {\mathcal K},M_1\} = g_1 M_1 + g_2 \mathcal{I}, \quad \{ {\mathcal K},M_2\} = g_3 M_2 + g_4 \mathcal{I}.
\lab{cJ_rel} \ee
The parameters $g_i, i=1, \dots, 4$ are assumed to be real numbers and it is further imposed that  $M_1, M_2$ are involution operators, i.e. that
\be
M_1^2=M_2^2=\mathcal{I}.
\lab{MM_inv} \ee

Our main goal will be to construct representations of this algebra to derive in the process the eigenvalues $\lambda_n$ of the operator ${\mathcal K}$ and the Verblunsky parameters $a_n$ as entries of the matrices ${\mathcal M}_1, \: {\mathcal M}_2$ (see \eqref{M1_def}, \eqref{M2_def}). We will also seek in the next section, to obtain the most general expression for the Dunkl type differential operator ${\mathcal K}$ by considering realizations in terms of operators in the variable $z$. Our observations will justify referring to the above construct as the {\it circle Jacobi algebra}.

First, let us make some remarks concerning the ``canonical" form of the circle Jacobi algebra \re{cJ_rel}. The affine transformation ${\mathcal K} \to \mu {\mathcal K} + \nu \mathcal{I}$ with arbitrary parameters $\mu,\nu$ allows us to reduce the number of independent parameters to two. We shall say that the algebra \re{cJ_rel} is {\it nondegenerate} if $g_3 \ne g_1$ and $g_2 \ne 0$. In this case it is  possible to bring the defining relations to the following ``canonical" form (that proves convenient in connection with Jacobi OPUC): 
\ba
&&\{\mathcal{K}, M_1\} = (\alpha+\beta+1) \left( M_1-\mathcal{I}\right),  \nonumber \\
&&\{\mathcal {K}, M_2\} = (2+\alpha+\beta) M_2 +(\alpha-\beta)\mathcal{I}
\lab{KM12} \ea
with $\alpha,\beta$, two independent parameters.

We shall now explain how to obtain the parameters $\lambda_n$ and $a_n$ by constructing representations of the circle Jacobi algebra. We start from the standard form of a pair of involution operators $M_1,M_2$ in terms of semi-infinite symmetric block-diagonal matrices \re{M1_def} - \re{M2_def}. The matrix $\mathcal K$ is taken to be diagonal with entries $\{\lambda_0, \lambda_1, \dots \}$.
From the algebraic relations \re{KM12} for the non-diagonal entries of $M_1$ and $M_2$, we find the equations 
\be
\left( 1-a_{2n}^2\right) \left(\lambda_{2n}+\lambda_{2n+1}-\alpha-\beta-1 \right)=0, \quad 
\left( 1-a_{2n+1}^2\right) \left(\lambda_{2n+1}+\lambda_{2n+2}-\alpha-\beta  \right)=0, \; n=0,1,2,\dots
\lab{eq_lambda} \ee
Demanding that the corresponding OPUC be nondegenerate, i.e. that $a_n^2 \ne 1$ for $n=0,1,2,\dots$, 
\re{eq_lambda} yields the following expressions for $\lambda_n$:
\be
\lambda_{2n} = \lambda_0 - n, \quad \lambda_{2n-1} = \alpha+\beta +n+1-\lambda_0
\lab{lam_sol} \ee
where $\lambda_0$ is still an undetermined parameter. Turning to the diagonal entries of the relations \re{KM12}, one first finds that
\be
\lambda_0 = 0
\lab{lam0} \ee
and then obtains for the Verblunsky parameters the explicit formulas:
\be
a_n = - \frac{\alpha + 1/2+ (-1)^{n+1} \left(  \beta +1/2 \right) }{n+\alpha+\beta+2}, \quad n=0,1,2,\dots
\lab{an_sol} \ee
These Verblunsky parameters correspond to those of the Jacobi OPUC (see, e.g. \cite{Badkov}, \cite{DVZ}) that are orthogonal on the unit circle \cite{Badkov}
\be
\int_{0}^{2\pi} w(\theta) \Phi_n\left(e^{i \theta}\right) \Phi_m \left(e^{-i \theta}\right) d \theta = 0, \; n \ne m,
\lab{ort_JOPUC} \ee
with respect to the positive weight function
\be
w(\theta) = \left( 1-\cos \theta \right)^{\alpha+1/2}  \left( 1+ \cos \theta \right)^{\beta+1/2}.
\lab{w_JOPUC} \ee
This representation theoretic derivation of the Verblunsky parameters of the Jacobi OPUC justifies a posteriori the name circle Jacobi algebra given to the algebra defined by \eqref{cJ_rel} and \eqref{MM_inv}.
Notice that the single moment OPUC are a special case of the Jacobi OPUC with $\alpha=1/2, \beta=-1/2$.

\section{CMV bispectrality of Jacobi OPUC}
\setcounter{equation}{0}
In this section we present the eigenvalue equation for the generic Jacobi OPUC with Verblunsky parameters \re{an_sol} . This calls for a realization of the circle Jacobi algebra on the space of Laurent polynomials $f(z)$. The generators $M_1$ and $M_2$ are then represented by the operators already introduced in Section 3:
\be
M_1 = R, \; M_2 = zR,
\lab{M12R} \ee 
where $R$ is the reflection operator:
\be
R f(z) =f(1/z).
\lab{Rff} \ee
In this realization, the operator ${\mathcal K}$ should be some differential operator whose eigenvalues are the $\lambda_n$ \re{lam_sol} that have been obtained in the preceding section. In the special case of the single moment OPUC this operator is given by \re{L_SM}. 
For the generic Jacobi OPUC, this operator ${\mathcal K}$ turns out to be the following first-order operator
\be
\mathcal{K}= z \partial_z  + \frac{z\left((\alpha+\beta+1)z +\alpha-\beta \right)}{1-z^2} \left(R-\mathcal{I} \right)
\lab{KJop} \ee
and we have:
\begin{pr}  
The CMV Laurent polynomials $\psi_n(z)$ corresponding to the Jacobi OPUC with Verblunsky parameters \re{an_sol} satisfy the eigenvalue equation
\be
\mathcal{K} \psi_n(z) = \lambda_n \psi_n(z),
\lab{EIG_JOPUC} \ee
where 
\be
\lambda_n = \left\{ -n/2,  \; n \; \mbox{even} \atop  (n+1)/2+\alpha+\beta+1 , \; n \; \mbox{odd} . \right.
\lab{l_JOPUC} \ee
\end{pr}
{\it Proof}. It was already shown by Szeg\H{o} \cite{Szego} that the pair of PRL $P_n(x), Q_n(x)$ corresponding to the Jacobi OPUC coincides with the ordinary Jacobi polynomials $P_n^{(\alpha, \beta)}(x/2)$ and $P_n^{(\alpha+1, \beta+1)}(x/2)$ which are orthogonal on the interval $[-2,2]$ with the weight functions $(2-x)^{\alpha}(2+x)^{\beta}$ and $(2-x)^{\alpha+1}(2+x)^{\beta+1}$ respectively.

The identity
\be
z \partial_z P_n(x(z)) = n (z-1/z) Q_{n-1}(x(z))
\lab{rel_PQ} \ee
follows from the well known relation between Jacobi polynomials with parameters $(\alpha,\beta)$ and $(\alpha+1,\beta+1)$ \cite{KLS}.
Moreover, the differential equation that the Jacobi polynomials $P_n(x(z))$ satisfy \cite{KLS} can be presented in the form
\be
z^2 \partial_z^2 P_n(x(z)) + z\frac{(\alpha+\beta+2)z^2+2(\alpha-\beta)z +\alpha+\beta }{z^2-1} \partial_z P_n(x(z)) = n(n+\alpha+\beta+1) P_n(x(z)).
\lab{DEP} \ee

Now use \re{rel_PQ} to replace $Q_{n-1}(x(z))$ in \re{psi_PQ} with $\partial_z P_n(x(z))$ and note that $\psi_n(z)$ can then be presented as a linear combination of $P_n(x(z))$ and $\partial_z P_n(x(z))$. Applying the operator $\mathcal K$ to the formulas for $\psi_n(z)$  resulting from this substitution in \re{psi_PQ}, one obtains ${\mathcal K} \psi_n(z)$ as a linear combination of $P_n(x(z)), \partial_z P_n(x(z))$ and $\partial_z^2 P_n(x(z))$. Eliminating $\partial_z^2 P_n(x(z))$ with the help of \re{DEP}, one finally arrives at \re{EIG_JOPUC}.

It is interesting to remark that Dunkl type differential equations similar to \re{EIG_JOPUC} were derived in \cite{KB} for the so-called ``nonsymmetric Jacobi polynomials" $E_n(z), \: n=0, \pm 1, \pm 2, \dots$ in \cite{KB}. We therefore see that these nonsymmetric polynomials coincide with the CMV Laurent polynomials $\psi_n(z)$ corresponding to the Jacobi OPUC. Precisely in the notation of in \cite{KB}, the correspondance goes as follows:
\be
\psi_{2n}(z) = E_n(z), n=0,1,2,\dots, \quad \psi_{2n-1}(z) = E_{-n}(z), \: n=1,2,\dots
\lab{psi_E} \ee
Moreover, as mentioned by the authors of \cite{KB},  the operator $\mathcal K$ \re{KJop} coincides with Cherednik's Dunkl differential operator for the root system $BC_1$ \cite{Cher}.
This establishes connections between Cherednik's Dunkl trigonometric operator, nonsymmetric Jacobi polynomials, and CMV bispectrality for circle analogs of the Jacobi polynomials. As expected, 
it can moreover be  directly verified that the operators $M_1,M_2, \mathcal K$ acting as specified above on the space of Laurent polynomials $f(z)$, do satisfy the circle Jacobi algebra relations.
Summing up:  

(i) the Jacobi OPUC are CMV-bispectral;

(ii) their bispectrality is encoded in the circle Jacobi algebra.\\
{\it Remark}. Notice that the Dunkl type differential operator  \re{KJop} simplifies to $z \partial_z$ when $\alpha=\beta=-1/2$. This corresponds to the ``free" OPUC case because it follows from \re{an_sol} that $a_n=0, \: n=0,1,2,\dots$. This is in accordance with the result obtained in \cite{GV} and already mentioned in Section 5 according to which the ``free" polynomials on the circle are the only OPUC which satisfy an eigenvalue equation of the form \eqref{W_psi} with $W$ a pure differential operator in $z$ which must be the operator $z\partial_z$ (or polynomials thereof).

\section{Embedding of the Jacobi algebra in the circle Jacobi algebra}
\setcounter{equation}{0}
In this section we establish the relation between the circle Jacobi algebra \re{KM12} and the ordinary quadratic Jacobi algebra $QJ(3)$ \cite{GIVZ} which encodes the bispectrality of the Jacobi polynomials on the real line. This algebra can be presented as having two generators $K_1$ and $K_2$ obeying the commutation relations \footnote{The presentation of $QJ(3)$ given here differs from the one in \cite{GIVZ}, the third generator $K_3 = [K_1,K_2]$ has been eliminated by the introduction of double commutators and trivial redefinitions of the generators have been performed.}:

\ba
&&[K_1,[K_1,K_2]] =  a K_1^2 + d K_1 + e_1 \mathcal{I} \nonumber \\
&&[K_2,[K_2,K_1]]= a \{K_1,K_2\} + c K_1 + d K_2 + e_2 \mathcal{I},
\lab{JA} \ea
where $a, c, d, e_1, e_2$ are the structure constants of the algebra and $[A,B]=AB-BA$ denotes the commutator. Note that when $a=0$, the quadratic nature of $QJ(3)$ disappears as it reduces to a Lie algebra isomorphic to $sl_2$.

The importance of the Jacobi algebra $QJ(3)$ lies in the observation that it describes a ``hidden symmetry" of the Gauss hypergeometric equation and of the Jacobi polynomials. Let $\pi_k(x)$ a polynomial of degree $k$ in $x$. One can set $K_1=x$ and $K_2=\pi_2(x) \partial_x^2 + \pi_1(x) \partial_x + \pi_0$ which is a hypergeometric operator, and check that those identifications of $K_1$ and $K_2$ satisfy the Jacobi algebra relations \eqref{JA}. From representations of the Jacobi algebra one can derive the basic properties of the Jacobi polynomials (which are eigenfunctions of the operator $K_2$) \cite{GLZ}. Moreover, the Jacobi algebra provides an algebraic description of ``exactly solvable" potentials in quantum mechanics (see \cite{GLZ} for details). 

Consider now the following elements in the circle Jacobi algebra:
\be
X= M_2 M_1 + M_1 M_2, \quad Y= {\mathcal K}^2 -\left(\alpha+\beta+1 \right){\mathcal K}. 
\lab{XY} \ee
It is easily verified that both $X$ and $Y$ commute with $M_1$:
\be
[X,M_1]=[Y,M_1]=0.
\lab{CM_XY_M1} \ee
Then by a direct computation, one can derive the following relations
\be
X^2Y +YX^2 -2XYX = [X,[X,Y]]=2 X^2 -8 \mathcal{I}
\lab{JR1} \ee 
and
\ba
&&Y^2X +XY^2 -2YXY =[Y,[Y,X]] = \nonumber \\
&&2 \{ X,Y\} + \left( \alpha+ \beta\right) \left( \alpha+ \beta +2 \right) X + 2(\beta-\alpha) M_1 +2\left( \alpha-\beta\right)  \left( \alpha+\beta+1\right)    \mathcal{I}.
\lab{JR2} \ea

 Comparing \re{JR1} and \re{JR2} with the defining relations of the quadratic Jacobi algebra \re{JA}, and identifying $X=K_1$ and $Y=K_2$, we observe that the operators $X,Y$ realize a {\it central extension} of the Jacobi algebra $QJ(3)$. Indeed we find that the right hand side of \eqref{JR2} contains with respect to the second relation \eqref{JA}, an extra operator $2(\beta-\alpha) M_1$ which commutes with both generators $X$ and $Y$ and can be incorporated with the identity operator in an irreducible representation. Such central extensions of quadratic algebras have arisen in various contexts one being in the comparison by Koornwinder \cite{Koo} of the Askey-Wilson algebra and the double affine Hecke algebras of rank one.

Having related the circle Jacobi algebra with the quadratic Jacobi algebra $QJ(3)$, there remains to consider the eigenvalue problem for the operator $Y$
\be
Y \psi_n = \Lambda_n \psi_n,
\lab{Y_psi} \ee
where
\be
\Lambda_n = \lambda_n^2 -(\alpha+\beta+1)\lambda_n.
\lab{L_l} \ee
Using the explicit expression \re{lam_sol} for $\lambda_n$, we have
\be
\Lambda_{2n} = n(\alpha+\beta+n+1), \; \Lambda_{2n-1} = \Lambda_{2n}.
\ee 
It then follows from formulas \re{PQ_psi} that
\be
Y P_n(x(z)) = \Lambda_{2n} P_n(x(z)), \quad Y F_n(z) = \Lambda_{2n} F_n(z),
\lab{Y_PQ} \ee  
where
\be
F_n(z) = \left(z-z^{-1} \right) Q_{n-1}(z).
\lab{F_n} \ee
We thus find that the Jacobi polynomials $P_n(x(z))$ are eigenfunctions of the differential operator $Y$ corresponding to the eigenvalue $R=1$ of the reflection operator $R$. This is obvious because the polynomials $P_n(x(z))$ are symmetric, i.e.
\be
RP_n(x(z)) = P_n\left(x\left(z^{-1}\right )\right) =P_n(x(z)).
\lab{RPP} \ee
The second eigensolution $F_n(z)$ corresponds to $R=-1$. This is also obvious because $F_n(x(z))$ is antisymmetric 
\be
RF_n(x(z)) = F_n\left(x\left(z^{-1}\right )\right) =-F_n(x(z)).
\lab{RFF} \ee   
Hence the symmetric and antisymmetric eigenfunctions $P_n(x(z))$ and $F_n(x(z))$ of the operator $Y$ have the same eigenvalue $\Lambda_{2n}$ and are distinguished by the eigenvalue $(\pm 1)$ of the reflection operator $R$.

\section{Conclusion}

Let us conclude by underscoring the two main results of this report.

\vspace{2mm}

1. We have shown that the Jacobi OPUC are ``CMV-classical". That is, they satisfy a 5-term recurrence relation of CMV type (like all OPUC) and are in addition eigenfunctions of a first order differential operator of Dunkl type.

\vspace{2mm}

2. We have presented an algebraic interpretation of the Jacobi OPUC. This involved the introduction of the circle Jacobi algebra that was shown to encode the bispectral properties of the Jacobi OPUC in a way very similar to how the ordinary Jacobi algebra encapsulates the features of the Jacobi polynomials on the real line. Indeed, it was seen that a complete characterization of the Jacobi OPUC could be obtained from the construction of representations of the circle Jacobi algebra.

 \vspace{2mm}

The Jacobi OPUC thus stand as a very instructive example of ``classical" polynomials on the unit circle. It is expected that there are many other such ``classical" OPUC. It might hence be appropriate to revisit the commonly held opinion \cite{Simon} according to which the OPUC are poorer than the PRL with respect to the ``classical" property.  An open problem is obviously to classify all the CMV-classical OPUC. One could speculate that there might exist a circle analog of the Askey scheme. This exploration is part of our future plans.

\vspace{23mm}

{\large \bf Acknowledgments.} 

AZ thanks the Centre de Recherches Math\'ematiques (Universit\'e de
Montr\'eal) for hospitality.   
The research of LV is supported in part by a research grant from the Natural Sciences and Engineering Research Council
(NSERC) of Canada.

\vspace{5mm}

\bb{99}

%\bi{Charris} B.H.Aldana, J.A.Charris and
%O.Mora-Valbuena, {\it  On Block Recursions,
%Askey's Sieved Jacobi Polynomials
%and two related Systems}, Colloqium Mathematicum
%{\bf 78} 1998, 57--91.

%\bi{AAA}  W.Al-Salam, W. R. Allaway and R.Askey, {\it Sieved ultraspherical polynomials}, Trans. Amer. Math. Soc. {\bf 284} (1984), %39--55.

\bi{Askey} R. Askey, {\it R. Askey, Orthogonal polynomials old and new, and some combinatorial connections}, in: Enumeration and Design, D. M. Jackson and S. A. Vanstone (eds.), Academic
Press, Toronto, Ont., 1984, 67--84.

\bi{Badkov} V. M. Badkov, {\it Systems of orthogonal polynomials explicitly represented by the Jacobi polynomials}, Mathematical
Notes of the Academy of Sciences of the USSR {\bf 42}, 858--863 (1987).

%\bi{BG} F.Bouzeffour and M.Garayev, {\it Fractional supersymmetric quantum mechanics and lacunary Hermite polynomials}, Analysis and %Math.Physics, {\bf 11} (1921), article 17. arXiv:1810.08275v2.

\bi{Cher}  I. Cherednik, {\it Double affine Hecke algebras, Knizhnik-Zamolodchikov equations, and Macdonald's operators}, Int. Math. Res. Not. (1992), no. 9, 171--180.

\bi{DGTVZ} H. De Bie, V. X. Genest, S. Tsujimoto, L. Vinet, A. Zhedanov, {\it The Bannai-Ito algebra and some applications}, J. Phys.: Conf. Ser, {\bf 597} (2015) 012001

\bi{DVZ} M. Derevyagin, L. Vinet, A. Zhedanov, {\it CMV matrices and Little and Big -1
Jacobi Polynomials}, Constr. Approx. {\bf 36} (2012), 513--535.

%\bi{DO} C.F.Dunkl and E.M.Opdam, {\it Dunkl operators for complex reflection groups}, Proc.London Math.Soc. , {\bf 86} , (2003) , %70--108. 
\bi{DG} J. J. Duistermaat, F. A. Gr\"unbaum, {\it Differential Equations in the Spectral Parameter}, Commun. Math. Phys. {\bf 103} (1986), 177-240

\bi{GIVZ} V. Genest, M. E. H. Ismail, L. Vinet, A. Zhedanov, {\it Tridiagonalization of the hypergeometric operator and the Racah--Wilson algebra}, Proc. Amer. Math. Soc. {\bf 144} (2016), 4441--4454.

\bi{Ger} Ya. L. Geronimus,\quad {\it Polynomials Orthogonal on a
Circle and their Applications}, \\ Am.Math.Transl.,Ser.1 {\bf
3}(1962), 1-78.

\bi{GLZ} Ya. I. Granovskii, I. M. Lutzenko, A. S. Zhedanov, {\it Mutual integrability, quadratic algebras, and dynamical symmetry}, Ann. Physics,  {\bf 217} (1992), no. 1--20,

\bi{GV} F. A. Gr\"unbaum, L. Vel\'asquez, {\it The CMV Bispectral Problem}, IMRN, {\bf 2017}, (2017),  5833--5860. Arxiv 1607.01962v2

\bi{GVZ} V. X. Genest, L. Vinet, A. Zhedanov, {\it The Racah algebra and superintegrable models}, J. Phys.: Conf. Ser. {\bf 512} (2014), 012011

\bi{Ismail} M. E. H. Ismail, {\it Classical and Quantum orthogonal
polynomials in one variable}. Encyclopedia of Mathematics and its
Applications (No. 98), Cambridge, 2005.

%\bi{Ismail_Li} M.E.H.Ismail and X.Li, {\it On sieved orthogonal polynomials. IX: Orthogonality on the unit circle}, Pacific J.Math. %{\bf 153} (1992), 289--297.

%\bi{HM}  A. Herrera-Poyatos, Pieter Moree, {\it Coefficients and higher order derivatives of cyclotomic polynomials: Old and new}, Expos. Math. %(to appear), ArXive 1805.05207.

\bi{KLS} R. Koekoek, P. A. Lesky, R. F. Swarttouw. Hypergeometric orthogonal polynomials
and their q-analogues. Springer Science \& Business Media, 2010.

\bi{Koo} T. Koornwinder, {\it The Relationship between Zhedanov's Algebra $AW(3)$
and the Double Affine Hecke Algebra in the Rank One Case}, SIGMA {\bf 3} (2007), 063.

\bi{KB} T. Koornwinder, F. Bouzeffour, {\it Nonsymmetric Askey-Wilson polynomials as vector-valued polynomials}, Appl. Anal. {\bf 90}, 2011, 731--746. ArXiv:1006.1140v3.

%\bi{Kron} L.~Kronecker, {\it Zwei S\"atze 
%\"uber Gleichungen mit ganzzahligen Coefficienten}, Crelle, Oeuvres
%I (1857) 105--108.

%\bi{Lang} S.Lang, Algebra.  

\bi{Simon} B. Simon, {\it Orthogonal Polynomials On The Unit
Circle}, AMS, 2005.

\bi{Simon2} B. Simon, {\it OPUC on one foot}, Bull.Am.Math.Soc., {\bf 42}, (2005), 431--460.

\bi{Szego} G. Szeg\H{o}, {\it Orthogonal Polynomials}, American Mathematical Society, 1939. 4th Edition, 1975

\bi{VZ_Bochner} L. Vinet, A. Zhedanov, {\it Generalized Bochner theorem: Characterization of the Askey-Wilson polynomials}, J. Comp. Appl. Math. {\bf 211}, (2008), 45--56.

\bi{VZ_Askey} L. Vinet, A. Zhedanov, {\it An algebraic treatment of the Askey biorthogonal polynomials on the unit circle}, 
Forum of Mathematics, Sigma , {\bf 9} , 2021 , e68. arXiv:2102.01779v1.

\bi{Zhe_AW}  A.S.Zhedanov, {\it ``Hidden symmetry'' of Askey-Wilson polynomials}, Theoret. and Math. Phys. 89 (1991),
1146--1157.

\eb

\end{document}